\newif\ifenglish
\englishtrue
\input AHTOHFIG.STY
\raggedbottom
\def\KL{{\bf KL}}
\def\KLs{{\bf KLs}}
\def\Mu{{\rm M}}
\def\Alpha{{\rm A}}

\count0=1
\headline{\ifnum\count0=1 \vtt UDC 512.543.7\hss\else\hss\fi}

\centerline{\ssdbf THE KERVAIRE--LAUDENBACH CONJECTURE
AND PRESENTATIONS OF SIMPLE GROUPS%
\footnote{*}{\rm This work was supported by the Russian Foundation for
Basic Research, project no. 02-01-00170.}}

\smallskip
\centerline{\ss Anton A. Klyachko}
\smallskip
{
\ssqi
\centerline{Faculty of Mechanics and Mathematics, Moscow State University}
\centerline{Moscow 119992, Leninskie gory, MSU}
\centerline{klyachko@daniil.math.msu.su}
}
\medskip
{\narrower\small\noindent
The statement ``no nonabelian simple group can be obtained from a nonsimple
group by adding one generator and one relator"
\item{1)} is equivalent to the Kervaire--Laudenbach conjecture;
\item{2)} becomes true under the additional assumption that the initial
nonsimple group is either finite or torsion-free.

\noindent{\it Key words\/\rm:} Kervaire--Laudenbach conjecture, relative
presentations, simple groups, car motion, cocar comotion.

\noindent{\it AMS MSC\/\rm:} 20E32, 20F05, 20F06.

}

\s 1. Conventions

For a group $G=\gp{X\ |\ R}$, the symbol $\~G$ denotes the group
$\gp{G,t\ |\ w=1}\:=\gp{X\cup\{t\}\ |\ R\cup\{w\}}$ obtained from $G$
by adding one generator and one (arbitrary) relator.

In this paper, the term ``simple group" means ``nonabelian simple group".

\s 2. Introduction

Recall an old well-known unproven group-theoretic conjecture with
topological origin.

\proclaim{Kervaire--Laudenbach conjecture (\KL)} {\rm (see, e.g.,
[LS77], [MKS66])}.
If a group $G$ is nontrivial, then the group $\~G$ is also nontrivial.

Let us state a similar conjecture.

\proclaim{``Simple" Kervaire--Laudenbach conjecture (\KLs)}.
If a group $G$ is nonsimple, then the group $\~G$ is also nonsimple.

These problems are more than similar.

\Proposition 1. The conjecture \KL\ is equivalent to the
conjecture \KLs.

The best-known result partially confirming the conjecture \KL\
is the following remarkable theorem.

\Theorem \rm (Gerstenhaber and Rothaus [GR62]). If a group $G$ is
finite and nontrivial, then the group $\~G$ is also nontrivial.%
\footnote{**}{\rm Actually, the Gerstenhaber--Rothaus theorem is more
general; we formulate only the most important special case (see [LS77]).}

From this result, it is easy to derive the analogious ``simple" theorem.

\Theorem 1. If a group $G$ is finite and nonsimple, then
the group $\~G$ is also nonsimple.

Consider another result partially confirming the conjecture \KL.

\Theorem \rm [K93]. If a group $G$ is torsion-free and
nontrivial, then the group $\~G$ is also nontrivial.

The analogious ``simple" theorem is as follows.

\Theorem 2. If a group $G$ is torsion-free and nonsimple,
then the group $\~G$ is also nonsimple.

The proof of Theorem 1 and Proposition 1 is very short. The main contents
of this paper are a proof of Theorem~2. Actually, we establish a stronger
fact.

\proclaim{Main theorem}. If a group $G$ is torsion-free, then the group
$\~G$ is simple if and only if $G$ is simple and the word $w$ is
conjugate in the free product $G*\gp{t}_\infty$ to a word of the form
$t^{\pm1}g$, where $g\in G$.

We do not know whether or not the similar strengthening of Theorem 1 is
true.

Note that the main theorem and the well-known fact that any torsion-free
group embeds into a simple torsion-free group immediately implies the
following result of Cohen and Rourke.

\Theorem \rm [CR01]. If a group $G$ is torsion-free, then
the natural mapping $G\to\~G$ is surjective if and only if
$w$ is conjugate in the free product $G*\gp{t}_\infty$ to a word
of the form $t^{\pm1}g$, where $g\in G$.

Our proof of the main theorem almost immediately divides, depending on
the word $w$, into two cases --- easy and difficult. It is amusing that the
difficult case corresponds precisely to the words of complexity 1 (in the
sense of Forester and Rourke [FoR03]); it is much easier to deal with
words of higher complexity.

Among other things, this paper contains a complete and self-contained
exposition of the method of comotions, which we apply to prove the main
theorem in the difficult case, while in the easy case our geometric
argument is reduced to a simple application of the car-motion lemma from
[K93].


\s 3. Proof of Proposition 1 and Theorem 1

Let a cyclically reduced form of the word $w$ be
$$
w\equiv g_1t^{\epsilon_1}\dots g_nt^{\epsilon_n},\quad\hbox{where }
\epsilon_i\in\{\pm1\},\ g_i\in G.
$$
(These denotations are assumed to be fixed throughout the paper.)

The following reformulation of the Kervaire--Laudenbach conjecture is
well-known.

\Proposition 2 \rm(folklore). The conjecture \KL\ is equivalent to the
following conjecture \KL$'$.

\Conjecture \KL$'$. If $\sum\epsilon_i=\pm1$, then the natural mapping
$G\to\~G$ is injective.

\Proof
Let us suppose that \KL$'$ is true and prove \KL.
If $G\ne\1$ and $\sum\epsilon_i=\pm1$, then $\~G$ is nontrivial, because
it contains a nontrivial subgroup isomorphic to $G$. If
$\sum\epsilon_i\ne\pm1$, then $\~G$ is nontrivial, because it admits an
epimorphism onto the nontrivial group $\Z/(\sum\epsilon_i)\Z$.

Now, let us suppose that \KL$'$ is false and disprove \KL.
Let $\sum\epsilon_i=\pm1$, and let $N\ne\1$ be the kernel of the natural
mapping $G\to\~G$. It is well known that any group $G$ embeds into a
simple group $H$ (see, e.g., [LS77]). The group $\~H\:=\gp{H,t\ |\ w=1}$
is obviously trivial:
$$
\~H\iso\~H/\nc N \iso \~H/\nc H \iso
\Z\!\Bigm/\!\!\left(\hbox{$\sum\epsilon_i$}\right)\Z \iso\1.
$$
Here (and throughout this paper), $\nc X$ denotes the normal subgroup
generated by a set $X$. Proposition 2 is proven.

\smallskip
\noindent{\bf Proof of Proposition 1.}
Suppose that \KL\ is false, i.e., for some nontrivial group $G$, the
group $\~G$ is trivial. Let $S$ be a simple group. The group $G\times S$
is nonsimple, but $\gp{G\times S, t\ |\  w=1}\iso S$ is simple.
Thus, \KLs\ is also false.

Now, suppose that \KLs\ is false, i.e., for some nonsimple group $G$,
the group $\~G$ is simple. First of all, note that $G$ is nonabelian and
$\sum\epsilon_i=\pm1$, because otherwise $\~G$ would not coincide with
its commutator subgroup (or would be trivial) and, therefore, would not be
simple.  If the natural mapping  $G\to\~G$ is not injective, then \KL\ is
false by Proposition 2. Suppose that the natural mapping $G\to\~G$ is
injective.  Let $N$ be a proper nontrivial normal subgroup of $G$. Then
the quotient group $\~G/\nc{N}$, on the one hand, is trivial (because $\~G$ is
simple, and $N\ne\{1\}$ in $\~G$ by virtue of the injectivity of the
mapping $G\to\~G$), and, on the other hand, is isomorphic to the group
$\gp{G/N,t\ |\ w'=1}$, where $w'$ is obtained from $w$ by reduction modulo
$N$. Thus, \KL\ is also false. This completes the proof of Proposition 1.

\smallskip
\noindent{\bf Proof of Theorem 1.}
Theorem 1 can be derived from the Gerstenhaber--Rothaus theorem exactly
as \KLs\ is derived from \KL; we need only to mention that a quotient
group of a finite group is finite and that any finite group embeds into a
simple finite group (e.g., in the alternating one).

Note that it is not so easy to prove Theorem 2, because a quotient group
of a torsion-free group may have torsion. The remaining part of this
paper is the proof of Theorem 2 (to be more precise, of the main theorem).


\s 4. Algebraic lemmata

{\bf Notation} which we use is mainly standard. Note only that if
$k\in\Z$, $x$ and $y$ are elements of a group, and $\phi$ is a
homomorphism from this group to another, then $x^y$,
$x^{ky}$, $x^{-y}$, $x^\phi$, $x^{k\phi}$, and $x^{-\phi}$ denote
$y^{-1}xy$, $y^{-1}x^ky$, $y^{-1}x^{-1}y$, $\phi(x)$, $\phi(x^k)$, and
$\phi(x^{-1})$, respectively.

\Lemma 1.
Let $A$ and $B$ be torsion-free groups, and let $u\in (A*B)\setminus A$.
Then $\gp{A,u}=A*\gp{u}_\infty$. If, in addition, $A$ is nontrivial and
$B$ is noncyclic, then there exists a word $v\in A*B$ such that
$\gp{A,u,v}=A*\gp{u}_\infty*\gp{v}_\infty$.

\Proof
Clearly, we can assume that the first and last letters of the reduced form
of $u$ lie in $B$.
For such $u$, the first assertion of the lemma is obvious. Let
us prove the second assertion. If $u\in B$, then we can take $v=a^b$,
where $b$ is an arbitrary element of $B\setminus\gp u$ and $a$ is any
element of $A\setminus\1$. If $u\notin B$, $u=b_1a_1\dots b_k$, then we
can take $v=a^b$, where $b\in B\setminus\{b_1^{\pm1},b_k^{\pm1}\}$ and
$a\in A\setminus\{1\}$.  Lemma 1 is proven.

The following lemma is a version of an algebraical trick from [K93], which
was several times used for studying equations over groups and similar
matter (see [KP95], [CG95], [CG00], [CR01], [FeR96], [FeR98], [FoR03]). A
geometrical interpretation of this trick can be found in [FoR03].

\Lemma 2.
Suppose that the group $G$ is torsion-free,
$$
\sum\epsilon_i=1, \hbox{\quad and\quad } n>1.
\eqno{(1)}
$$
Then the group $\~G$ has a \(relative\)
presentation of the form
$$
\~G\iso\gp{H, t\ \Biggm|\ \{p^t=p^\phi,\ p\in P\setminus\1\},\
ct\prod_{i=0}^m(b_ia_i^t)=1},
\eqno{(2)}
$$
where $H$ is a group, $a_i,b_i,c\in H$, $P$ and $P^\phi$ are isomorphic
subgroups of $H$, $\phi\:P\to P^\phi$ is an isomorphism,
and the following conditions hold:
\item{\rm1)} $m\ge 0$ \(i.e., the product in $(2)$ is nonempty\);
\item{\rm2)} $a_i\notin P$ and
             $b_i\notin P^\phi$;
\item{\rm3)} $\gp{P,a_i}=P*\gp{a_i}_\infty$ and
             $\gp{P^\phi,b_i}=P^\phi*\gp{b_i}_\infty$ in $H$;
\item{\rm4)} if $G$ is noncyclic and $P\ne\1$, then, for each
           $i$, there exist elements $a'_i,b'_i\in H$ such that
           $\gp{P,a_i,a'_i}=P*\gp{a_i}_\infty*\gp{a'_i}_\infty$ and
           $\gp{P^\phi,b_i,b'_i}=P*\gp{b_i}_\infty*\gp{b'_i}_\infty$ in $H$;
\item{\rm5)} the groups $H$, $P$, and $P^\phi$ are free products of
             finitely many isomorphic copies of $G$:
             $H=G^{(0)}*\dots*G^{(s)}$, $P=G^{(0)}*\dots*G^{(s-1)}$,
             $P^\phi=G^{(1)}*\dots*G^{(s)}$, where $s\ge0$
             \(for $s=0$, $P$ and $P^\phi$ are assumed to be trivial\/\),
             and the isomorphism $\phi$ is the shift:
             $\left(G^{(i)}\right)^\phi=G^{(i+1)}$.

\Proof
First, let us show that the group $\~G$ has a presentation of the form (2)
satisfying condition 5. Since $\sum\epsilon_i=1$, the word $w$ can be
written in the form
$$
w=\left(\prod g_i^{t^{k_i}}\right) t.
$$
Conjugating, if necessary, $w$ by $t$, we can assume that $k_i\ge0$.
We set $g^{(i)}=g^{t^i}$ for $g\in G$, $G^{(i)}=G^{t^i}$,
$s=\max k_i$, and $c=\prod g_i^{(k_i)}$. We see that $\~G$ has
the presentation
$$
\~G\iso\gp{G^{(0)}*\dots*G^{(s)}, t\ \Biggm|\
\left\{\left(g^{(i)}\right)^t=g^{(i+1)},\ i=0,\dots,s-1,\ g\in G\right\},\
ct=1},
$$
which is a presentation of the form (2) (with $m=-1$) satisfying condition
5.

Now, from all presentations of the form (2) satisfying condition 5,
we choose presentations with minimal $s$, and from all these presentations
with minimal $s$, we choose one with minimal $m$. The obtained
presentation (2) is as required.

Indeed, if $m<0$ (i.e., $w=ct$, where $c\in H$), then $s=0$,
because otherwise we might decrease $s$ replacing
all fragments $g^{(s)}$ in the word $c$ by $(g^{(s-1)})^t$. But the
conditions $m<0$ and $s=0$ mean that the initial word $w$ has the form
$w=ct$, where $c\in G$, which contradicts the assumption $n>1$. Thus,
condition 1 holds.

Condition 2 holds because otherwise in presentation (2)
we might replace a fragment $t^{-1}a_it$ with $a_i\in P$
(or a fragment $tb_it^{-1}$ whith $b_i\in P^\phi$) by $a_i^\phi$
(or by $b_i^{\phi^{-1}}$, respectively), thereby decreasing $m$ (and not
increasing $s$).

Conditions 3 and 4 follow from conditions 2 and 5 by Lemma 1.
Lemma 2 is proven.


\s 5. Maps and motions

Throughout this paper, the term ``surface" means ``closed oriented
two-dimensional surface".

A {\it map} $\Mu$ on a surface $S$ is a finite set of continuous mappings
$\{\mu_i\:D_i\to S\}$, where $D_i$ is a compact oriented two-dimensional
disk, called the $i$th {\it face}, or {\it cell}, of the map; the boundary
of each face $D_i$ is partitioned into finitely many intervals
$e_{ij}\subset \d D_i$, called the {\it pre-edges} of the map, by a nonempty
set of points $c_{ij}\in \d D_i$, called the {\it corners} of the map.
The images of the corners $\mu_i(c_{ij})$ and pre-edges $\mu_i(e_{ij})$
are called the {\it vertices} and {\it edges} of the map, respectively.
It is assumed that
\item{1)}
 the restriction of $\mu_i$ to the interior of each
 face $D_i$ is a homeomorphic embedding preserving orientation; the
 restriction of $\mu_i$ to each pre-edge is a homeomorphic embedding;
\item{2)}
 different edges do not intersect;
\item{3)}
  the images of the interiors of different faces do not intersect;
\item{4)}
 $\bigcup\mu_i(D_i)=S$.

\noindent
Sometimes, we interpret a map $\Mu$ as a continuous mapping
$\Mu\:\coprod D_i\to S$ from a discrete union of disks onto the surface.

The union of all vertices and edges of a map is a graph on the surface,
called the {\it $1$-skeleton}.
The {\it multiplicity} of a point of the 1-skeleton is the number of edges
incindent to this point if this point is a vertex; if this point lies on an
edge, then its multiplicity is assumed to be two. In other words, the
multiplicity of a point $p$ is $\card{\Mu^{-1}(p)}$.

We say that a corner $c$ is a corner at a vertex $v$ if $\Mu(c)=v$.  There
is a natural cyclic order on the set of all corners at a vertex $v$; we
call two corners at $v$ {\it adjacent} if they are neighboring with
respect to this order.

By abuse of language, we say that a point or a subset of the surface is
contained in a face $D_i$ if it lies in the image of $\mu_i$. Similarly,
we say that a face $D_i$ is contained in some subset $X\subseteq S$ of
the surface $S$ if $\Mu(D_i)\subseteq X$.

Figure 1 presents a map on the sphere with 5 faces --- $A$, $B$, $C$,
$D$, and $E$, 18 corners --- $a_i$, $b_i$, $c_i$, $d_i$, and $e_i$, 6
vertices, 9 edges, and 18 pre-edges. Note that the number of corners always
equals the number of pre-edges and is twice the number of edges,
and the value
$$
e(S)\:=(the\ number\ of\ vertices)-(the\ number\ of\ edges)+
(the\ number\ of\ faces)
$$
does not depend on the choice of a map on the surface $S$ and is called
{\it the Euler characteristic} of this surface. The Euler characteristic
of the sphere (the only surface of our real interest in this paper)
is two.

\input 1.PIC

{\it A motion} on a surface $S$ is a map $\Mu\:\coprod D_i\to S$ on this
surface and a set of continuous mappings $\alpha_i\:\R\to\d D_i$.
The mapping $\alpha_i$ is called a {\it car} moving around the face $D_i$.
We say that a car $\alpha_i$ {\it is at} a corner $c\in\d D_i$ at a moment
of time $t\in\R$ if $\alpha_i(t)=c$; We also say that a car
$\alpha_i$ {\it is at} a point $p\in S$ at a moment $t\in\R$ if
$\mu_i(\alpha_i(t))=p$. If the number of cars being at moment $t\in\R$ at
a point $p$ of the 1-skeleton of $S$ equals the multiplicity of this
point (in other words, $\bigcup\alpha_i(t)\supseteq\Mu^{-1}(p)$), then we
say that at the point $p$ at the moment $t$ a {\it complete collision}
occurs; The point $p$ is called a {\it point of complete collision}.
Points of complete collision lying on edges are called simply {\it points
of collision}.

A motion is called {\it regular} if the mappings $\alpha_i$ are
orientation-preserving coverings. Simply speaking, in a regular
motion, a car moves around the boundary of its face anticlockwise
(the interior of the face remains on the left from the car)
without U-turns, stops, and ``infinite decelerations
and accelerations".

\Example 1.
On the map shown on Figure 1, we specify a regular motion as follows:
the cars $\alpha$, $\beta$, $\gamma$, $\delta$, and $\epsilon$
moving around the faces $A$, $B$, $C$, $D$, and $E$,
respectively, move with unit speed (one edge per unit time) in
the positive direction and, at the zero moment of time, the cars visit the
corners $a_0$, $b_0$, $c_0$, $d_0$ and $e_0$, respectively (at corners
with number $i$, the cars are at the moment $t=i$). In Figure 1, the
positions of the cars and the direction of their motion at the moment
$t=4/3$ are shown (the car $\beta'$ should be ignored for a while). This
motion has 3 points of complete collision, they are marked by the
exclamation signs in Figure 1: at the moments beeing multiples of 3,
complete collisions of the cars $\gamma$, $\delta$, and $\epsilon$ occur; at
the moments $t\in 3/2+6\Z$, the cars $\beta$ and $\epsilon$ collide on an
edge; in addition, at the moments beeing multiples of 3, the car $\alpha$
is at the dead end, this is also a complete collision, according to the
definition. We can vary slightly the schedule of the motion and reduce the
number of complete collision points to two (e.g., we can make the car
$\epsilon$ to move with speed 2 on the edges $[e_0,e_1]$ and $[e_1,e_2]$
and with speed 1/2 on the edge $[e_2,e_0]$). The further optimization of
the schedule is impossible, as the following lemma shows.

\Lemma 3 {\rm[K93] (see also [FeR96])}. Any regular motion on the
sphere has at least two points of complete collision.

Sometimes (see [K93] or [FeR96]), it is useful to consider motions
slightly more general than regular.

We call a continuous mapping $\alpha\:X\to Y$ from an oriented line or
circle $X$ to an oriented circle $Y$ (locally) {\it nondecreasing} if
the preimage of any interval $U\subset Y$ is a union of such intervals
that the restriction of $\alpha$ to each of these intervals is
a nondecreasing function (in the usual sense, as a function from one
oriented interval to another). We call a mapping $\alpha\:X\to Y$
from a line to a circle {\it proper} if the image of any
half-line $U\subset X$ is the entire circle $Y$.

The preimage of a point under a proper nondecreasing mapping from a line
to a circle is a discrete union of points and closed intervals.
A point $y\in Y$ whose preimage is nondiscrete is called a
{\it stop point} of the mapping.

A motion on a surface $S$ is called a {\it motion with separated
stops} if every car is a proper nondecreasing function each of whose stop
points is a corner (i.e., simply speaking, each car moves without U-turns
and infinite decelerations and accelerations moving around the
boundary of its face anticlockwise, possibly stopping for a finite time at
some corners); and there exists a set of corners called the
{\it stop corners} such that
\item{1)} the cars stop only at stop corners
          (possibly, at some stop corners, the cars do not stop);
\item{2)} at each vertex $v$ having stop corners at it,
          the stops are separated in the following sense: let
          $c_1,\dots,c_k$ be all stop corners at $v$ enumerated
          anticlockwise; it is required that, for each $i$, at corners
          $c_{i}$ and $c_{i+1}$ (subscripts are modulo $k$), cars are
          never located simultaneously. (In particular, this implies that
          $k\ge 2$.)

\Lemma 4 \rm{[K93] (see also [FeR96])}. Any motion with
separated stops on the sphere has at least two points of complete
collision.

\Proof
First, note that the motion of a car moving via a stop corner
can be slightly changed in a small neighborhood of this corner in such a
way that this car does stop at this corner, no new points of complete
collision arise, and the separated stops condition is not violated.

Now, assuming that each car does stop every time moving via a stop
corner, we make the following transformation of the map
called {\it blowing-up of stop corners} (see Fig. 2): we take a vertex
$v$ at which there are $k\ge2$ stop corners $c_1,\dots,c_k$; for each $i$,
we cut the surface along a small arc drawn from the vertex $v$ ``inside"
the corner $c_i$; the boundary of this cut is denoted $x_i$ (the left
boundary if to look from $v$) and $y_i$ (the right boundary if to look from
$v$). When we have made these cuts for all $i=1,\dots,k$, there
appears a hole on the surface with boundary consisting of the consecutive
segments $y_1,x_1,y_2,x_2,\dots,y_k,x_k$; let us remove this hole by glueing
each segment $x_i$ with $y_{i+1}$ (the subscripts are modulo $k$); we
obtain a new map $\Mu'$ on the same surface. This new map has the
additional pre-edges $x_i$ and $y_i$ (instead of the stop corners $c_i$)
and the additional edges $\Mu'(x_i)=\Mu'(y_{i+1})$.

\input 2.PIC

Having made the transformation described above for each vertex $v$ with
stop corners at it, we specify a regular motion on the obtained map as
follows: the cars move as on the initial map, but, during the time when a
car on the initial map is staying at a stop corner $c_i$, the
corresponding car on the new map moves uniformly along the pre-edges $x_i$
and $y_i$.  The stops separated condition implies that this regular
motion on the new map has no complete collisions on the additional edges
and at their ends.  Thus, the assertion of the lemma follows from Lemma 3.

Note that the separated stops condition means, in particular, that
no complete collision can occur at a vertex at which there is at least one
stop corner.


\s 6. Howie diagrams

Suppose that we have a map $\Mu$ on a surface $S$, the corners of the
map are labeled by elements of a group $H$, and the edges are oriented (in
the figures, we draw arrows on the edges) and labelled by elements of a
set $\{t_1,t_2,\dots\}$ disjoint from the group $H$. The label of a corner
or an edge $x$ is denoted by $\lambda(x)$.

The {\it label of a vertex} $v$ of such a map is defined by the formula
$$
\lambda(v)=\prod_{i=1}^k \lambda(c_i),
$$
where $c_1,\dots,c_k$ are all corners at $v$ listed clockwise.
The label of a vertex is an element of the group $H$ determined up to
conjugacy.

The {\it label of a face} $D$ is defined by the formula
$$
\lambda(D)=\prod_{i=1}^k
\bigl(\lambda(\Mu(e_i))\bigr)^{\epsilon_i}\lambda(c_i),
$$
where $e_1,\dots,e_k$ and $c_1,\dots,c_k$ are all pre-edges and all
corners of $D$ listed anticlockwise, the endpoints of $e_i$ are
$c_{i-1}$ and $c_i$ (subscripts are
modulo $k$), and $\epsilon_i=\pm1$ depending on whether the homeomorphism
$e_i\mathop\to\limits^\Mu\Mu(e_i)$ preserves or reverses orientation.
Simply speaking, to obtain the label of a face, we should go around its
boundary anticlockwise, writing out the labels of all corners and edges we
meet, the label of an edge traversed against the arrow should be raised to
the power $-1$.

The label of a face is an element of the group $H*F(t_1,t_2,\dots)$ (the
free product of $H$ and the free group with basis $\{t_1,t_2,\dots\}$)
determined up to a cyclic permutation. More precisely, the right-hand side
of our formula for $\lambda(D)$ is called the {\it label of the face $D$
written starting with the pre-edge $e_1$}.

Such a labelled map is called a {\it Howie diagram} (or
simply {\it diagram}) over a relative presentation
$$
\gp{H,t_1,t_2,\dots\ |\ w_1=1,w_2=1,\dots}
\eqno{(*)}
$$
if
\item{1)}
  some vertices and faces are separated out and called {\it
  exterior}, the remaining vertices and faces are called {\it interior};
\item{2)}
    the label of each interior face is a cyclic permutation of one of
  the words $w_i^{\pm1}$;
\item{3)}
  the label of each interior vertex is the identity element of $H$.

A diagram is said to be {\it reduced} if it contains no such
edge $e$ that both faces containing $e$ are interior, these faces
are different and their labels written starting with the $\Mu$-preimages of
$e$ are mutually inverse; such a pair of faces with a common edge is
called a {\it reducible pair}.

The following lemma is an analog of the van Kampen lemma for relative
presentations.

\Lemma 5 {\rm[H83]}.
The natural mapping from a group $H$ to the group with relative
presentation $(*)$ is noninjective if and only if there exists a
spherical diagram over this presentation with no exterior faces and a
single exterior vertex whose label is not 1 in $G$. A
minimal \(with respect to the number of faces\) such diagram is
reduced. \hfil\break
If this natural mapping is injective, then we have
the equivalence: the image of an element $u\in
H*F(t_1,t_2,\dots)\setminus \1$ is 1 in the group $(*)$ if and only
if there exists a spherical reduced diagram over this presentation without
exterior vertices and with a single exterior face with label $u$. A
minimal \(with respect to the number of faces\) such diagram is
also reduced.

Diagrams on the sphere with a single exterior face and no exterior
vertices are also called {\it disk diagrams}, the boundary of the exterior
face of such a diagram is called the {\it contour} of the diagram.

Let $\phi\:P\to P^\phi$ be an isomorphism between two subgroups of a group
$H$.  A relative presentation of the form
$$
\gp{H,t\ |\ \{p^t=p^\phi;\
p\in P\setminus\1\}, w_1=1,\ w_2=1,\ \dots}
\eqno{(**)}
$$
is called a {\it
$\phi$-presentation}. A diagram over a $\phi$-presentation $(**)$ is
called {\it $\phi$-reduced} if it is reduced and different interior cells
with labels of the form $p^tp^{-\phi}$, where $p\in P$, have no common
edges.

\Lemma 6. A minimal \(with respect to the number of faces\) diagram
among all spherical diagrams over a given $\phi$-presentation without
exterior faces and with a single exterior vertex with nontrivial label is
$\phi$-reduced. If no such diagrams exists, then a minimal diagram
among all disk diagrams with a given label of contour is $\phi$-reduced.
\rm In other words, the complete $\phi$-analog of Lemma 5 is valid.

\Proof
Indeed, if in some diagram over presentation $(**)$ a pair of cells
corresponding to relations of the form $p^tp^{-\phi}$, where $p\in P$ has
a common edge, then either this pair of cells is a reducible pair or
we can remove this common edge (Fig. 3) and obtain a diagram
with a smaller number of cells and with the same labels of exterior faces
and vertices; this means that the initial diagram is not minimal and
proves the lemma.

\input 3.PIC

A relative presentation ($\phi$-presentation) over which there exists
no reduced (respectively, $\phi$-reduced) spherical diagrams with a single
exterior vertex are called {\it aspherical} (respectively, {\it
$\phi$-aspherical}).

\Lemma 7.
Suppose that the group $H$ in presentation $(*)$
is nontrivial and there exists a word $w\in H*F(t_1,t_2,\dots)$
nonconjugate to the elements of $H\cup \{w_i^{\pm1}\}$ in
$H*F(t_1,t_2,\dots)$ and such that the presentation
$$
L=\gp{H,t_1,t_2,\dots\ |\ w=1,w_1=1,w_2=1,\dots}
$$
is aspherical (or $\phi$-aspherical, if
the initial presentation $(*)$ is a $\phi$-presentation).
Then the group $K$ with presentation $(*)$ is nonsimple.

\Proof
The ($\phi$-)asphericity of the presentation $L$ implies
the ($\phi$-)asphericity of presentation $(*)$.

Let us show that $w\ne1$ in the group $K$. Indeed, otherwise, by virtue of
the ($\phi$-)asphericity of $(*)$ and Lemmata 5 and 6, over this presentation
there would exist a ($\phi$-)reduced disk diagram whose contour label is
$w$.  But such a diagram can be considered as a ($\phi$-)reduced spherical
diagram over the presentation $L$ without exterior faces and vertices,
which contradicts the asphericity of $L$.

The natural mapping from $H$ to $L$ is injective by Lemma 5 (and
Lemma 6). Thus, $L=K/\nc w\ne\1$, and $\nc w$ is a proper nontrivial normal
subgroup of $K$. Q.E.D.


\s 7. Standard motion

Consider a map whose edges are oriented
(e.g., a Howie diagram). Such a map can have corners of four kinds
$(++)$, $(--)$, $(+-)$, and $(-+)$ (Fig. 4).
\noindent
\input 4.PIC

The following lemma is obvious.

\Lemma 8. In the anticlockwise listing of the corners at a vertex
$v$, the corners of type $(++)$ alternate with corners of type $(--)$. If
at a vertex $v$ there are no corners of type $(++)$, or, equivalently,
there are no corners of type $(--)$, then either all corners at $v$ are of
type $(+-)$ \(in this case, $v$ is called a {\it sink}\), or all corners at
$v$ are of type $(-+)$ \(in this case, $v$ is called a {\it source}\) \rm
(Fig. 5).
\noindent
\input 5.PIC

We say that a map with oriented edges is of type $A_m$, $m\ge0$, if the
sequence of pre-edge orientations of each face has one of the following
four forms:
\item{a)} $+-$ (Figure 6a);
\item{b)} $+(+-)^{m+1}$ (Figure 6b);
\item{c)} $-(-+)^{m+1}$ (Figure 6c);
\item{d)} $(+)^{k+1}(-)^{l+1}$, where $k,l\ge1$ (Figure 6d).

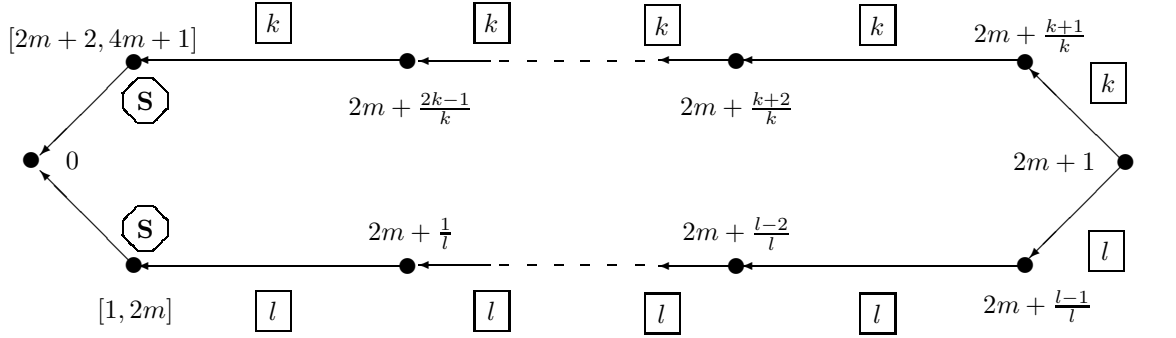
\input 6a.PIC
\input 6b.PIC
\input 6c.PIC
\input 6d.PIC

\noindent
We define the {\it standard motion} on a map of type $A_m$ as follows:

\item{a)}
  the car going around a face of type $+-$ moves anticlockwise
  uniformly with unit speed (one edge per unit time)
  visiting the corner of type $(+-)$ at the even moments of time
  (Fig. 6a);

\item{b)}
  for $m>0$, the car moving around a face of type $+(+-)^{m+1}$ stays at
 the corner of type $(++)$ during the time intervals
 $[2m+2,4m+1]+(4m+2)\Z$, and moves anticlockwise uniformly with unit speed
 all the remaining time; for $m=0$, such a car moves without stops with
 speed 2 on the positive pre-edges and with speed 1 on the negative ones
 visiting the corner of type $(+-)$ at the even moments of time (Fig. 6b);

\item{c)}
  for $m>0$, the car moving around a face of type $-(-+)^{m+1}$ stays at
 the corner of type $(--)$ during the time intervals $[1,2m]+(4m+2)\Z$,
 and go anticlockwise uniformly with unit speed all the remaining time;
 for $m=0$, such a car moves without stops with speed 2 on the negative
 pre-edges and with speed 1 on the positive ones visiting the corner of
 type $(+-)$ at the even moments of time (Fig. 6c);

\item{d)}
  for $m>0$, the car moving around a face of type
 $(+)^{k+1}(-)^{l+1}$ is at the corner of type $(+-)$ at time zero;
 then, it moves along the first negative pre-edge with unit speed;
 after that, it stops and waits during the time interval $[1,2m]$;
 next, it goes with speed $l$ along the remaining $l$ negative pre-edges
 and with speed $k$ along the $k$ positive pre-edges;
 after that, it stops and waits during the time interval $[2m+2,4m+1]$;
 then, it moves with unit speed through the last positive pre-edge and
 returns at the moment $4m+2$ again to the corner of type $(+-)$;
 after that, everything is repeated with period $4m+2$.
 For $m=0$ such a  car moves without stops with speed $l+1$ on the negative
 pre-edges and with speed $k+1$ on the positive ones visiting the
 corner of type $(+-)$ at the even moments of time (Fig. 6d).

\noindent
The standard motion is periodic with period $4m+2$ (on the faces of
type $+-$, the minimal period is 2). Figure 6 shows the detailed
schedule of the motion during the time interval $[0,4m+2)$, the
boxed numbers near edges indicate the speed of the car on these edges (by
default, the speed is 1). At corners labelled by the letter $\bf S$ the
car stops for time $2m-1$.

\Lemma 9. The standard motion on a map of type $A_m$ is a
motion with separated stops. The complete collisions can occur only
at vertices being sources or sinks and only at integer moments of time.

\Proof
Let us declare all corners of types $(++)$ and $(--)$ to be the stop
corners. Then, the stops are separated by virtue of Lemma 8 and the
fact that the schedule of the standard motion is such that cars are never
located simultaneously at corners of types $(++)$ and $(--)$ (the
corners of type $(--)$ are visited only during the first half of
the period, i.e, at moments from the intervals $(0,2m+1)+(4m+2)\Z$, while
the corners of type $(++)$ are visited during the second half of the period,
i.e., at moments from the intervals $(2m+1,4m+2)+(4m+2)\Z$).

A collision on an edge at a moment $t$ means that at this moment
the direction of the motion of one of the cars coincides with the
direction of the edge, while the direction of the motion of the other
colliding car is opposite to the direction of the edge (Fig. 7).
\nobreak
\input 7.PIC
\noindent
But the schedule of the standard motion is such that, at each moment $t$,
either all cars being on edges move in the direction of the edge (this is
so when the integer part of $t$ is odd), or all cars being on edges move
in the direction opposite to the direction of the edge (this is so when
the integer part of $t$ is even). Therefore, collisions can occur only at
vertices; the separatedness of stops implies that a vertex of complete
collision can not have stop corners and, therefore, is a source or a sink.
The cars visit such vertices only at integer moments of time (even for
sinks and odd for sources). The lemma is proven.


\s 8. Proof of the main theorem. The easy case

The ``if" part of the main theorem is obvious. Let us prove the
``only if" part. Clearly, the simplicity of $\~G$ implies that the group
$G$ coincides with its commutator subgroup and $\sum\epsilon_i=\pm1$
(because otherwise $\~G$ would not coincide with its commutator subgroup).
We suppose that $G$ coincides with its commutator subgroup and
conditions (1) hold; we have to prove that $\~G$ is nonsimple.

By Lemma 2, $\~G$ has presentation (2). The easy case considered in this
section is the case when $P\ne\1$ in presentation (2) or,
equivalently, the word $w$ is not conjugate to a word of the form
$ct\prod_{i=0}^m(b_ia_i^t)$, where $c, a_i, b_i\in G$. The nonsimplicity
of $\~G$ in this case is a corollary of Lemma 7 and the following fact.

\Lemma 10. Suppose that $G$ is a noncyclic torsion-free group, conditions
$(1)$ hold, and $P\ne\1$ in presentation $(2)$, then there exist elements
$a,b\in H$ such that the presentation
$$
\~G\!\!\Bigm/\!\!\nc{a^{t^2}b}\iso \gp{H, t\ \Biggm|\
\{p^t=p^\phi,\ p\in P\setminus\1\},\
ct\prod_{i=0}^m(b_ia_i^t)=1,\ a^{t^2}b=1 }
\eqno{(3)}
$$
obtained from presentation $(2)$ by adding the relator
$a^{t^2}b=1$ is $\phi$-aspherical.

\Proof
Let $a$ and $b$ be such elements of $H$ that
$\gp{P,a_m,a}=P*\gp{a_m}_\infty*\gp{a}_\infty$ and
$\gp{P^\phi,b_0,b}=P*\gp{b_0}_\infty*\gp{b}_\infty$.
Such elements exist by Lemma 2 (property 4).

We have to show that there is no $\phi$-reduced spherical diagram over
presentation (3) without exterior faces and with a single exterior vertex.
The cells of such a diagram have the forms shown in Figure 8. We see that
the diagram is a map of type $A_m$. Let us show that all complete
collisions of the standard motion on this map occur only at the exterior
vertex.

\input 8.PIC

According to Lemma 9, a complete collision can occur only at vertices
being sinks or sources.

Suppose that a vertex of complete collision is a sink. Then, all corners at
this vertex are of type $(+-)$; the label of each of these corners
is either $p^\phi$ (where $p\in P$), $b_i^{\pm1}$, or $b^{\pm1}$
(see Fig. 8). If at this vertex there are a corner labelled by $b_i^{\pm1}$
and a corner labelled by $b_j^{\pm1}$, where $i\ne j$, then a complete
collision does not occur at this vertex, because these two corners
are never visited simultaneously (compare Figs. 8, 6b, and 6c). If at
this vertex there is a corner labelled by $b^{\pm1}$, then a complete
collision can occur only at the moment $0$ ($\mod 4m+2$), because only in
this moment a car visit such a corner (Figs. 8 and 6d); but then, this
vertex of complete collision can not have corners with labels
$b_i^{\pm1}$, where $i\ne0$, because such corners are not visited at the
moment $0$ ($\mod 4m+2$) (Figs. 8, 6b, and 6c). Thus, the label of a
vertex being a sink at which a complete collision occurs is either
$$
\prod_j(b_i^{\epsilon_j}p_j^\phi)\quad\hbox{or}\quad
\prod_j(x^{\epsilon_j}p_j^\phi),
$$
where $\epsilon_j\in\Z$, $p_j\in P$, and $x\in\{b,b_0\}$.
But the label of an interior vertex must be 1. This means that we have
a nontrivial (because the diagram is $\phi$-reduced) relation of
the specified form, which contradicts property 3 from Lemma~2 and the
choice of $b$.

Similarly, supposing that a complete collision at a moment $t$ occurs at
an interior vertex being a source, we obtain a nontrivial relation of
the form
$$
\prod_j(a_i^{\epsilon_j}p_j)=1\ \ \hbox{(if $t\ne 2m+1\ (\mod 4m+2$))}
\quad\hbox{or}\quad
\prod_j(x^{\epsilon_j}p_j)=1\ \ \hbox{(if $t=2m+1\ (\mod 4m+2$))},
$$
where $\epsilon_j\in\Z$, $p_j\in P$, and $x\in\{a,a_m\}$,
which contradicts property 3 from Lemma 2 and the choice of $a$.

Thus, a complete collision can occur only at the exterior vertex,
but Lemma 4 says that there must be at least two different points of
complete collision. This contradiction completes the proof of Lemma 10 and
of the easy case of the main theorem.

\Remark. Note that as a byproduct we have proven the
following fact. \sl Suppose that we have a motion on a $\phi$-reduced
diagram over presentation $(2)$ satisfying properties $1$ and $3$ from
Lemma $2$ which is standard on the interior faces. Then complete
collisions can occur only at exterior vertices and on the boundaries of
exterior faces.
\rm In particular, this implies (by virtue of Lemmata 2,
4, 5, and 6) the main result of [K93]:
\sl if a group $G$ is
torsion-free and $\sum\epsilon_i=1$, then the natural mapping $G\to\~G$ is
injective.


\s 9. Comotions and multiple motions

In this section, we expose the theory of comotions, which were introduced
and applied to study equations over group in [K94] and [K97]. Our
present exposition differs from that of the works mentioned above
in terminilogy and generality.

The notion of a comotion is dual to that of a periodic motion.

\Def A {\it comotion} $\Alpha$ on a surface $S$ is a map
$\Mu\:\coprod D_i\to S$ on this surface and a set of continuous mappings
$\{\alpha_i\:\d D_i\to\bf T\}$, where $\bf T$ is an oriented circle,
called the circle of time. Sometimes, we shall interpret a comotion as
a continuous mapping $\Alpha\:\coprod\d D_i\to\bf T$. The mapping $\alpha_i$
is called the {\it cocar} moving around the face $D_i$. We say that a
cocar $\alpha_i$ {\it is at} a point $p\in S$ at time
$t\in\bf T$ if $\alpha_i(\mu_i^{-1}(p))\ni t$. We say that
a {\it complete collision} occurs at a point $p$ of the 1-skeleton of the
surface $S$ at a moment $t\in\bf T$ if
$\Alpha^{-1}(t)\supseteq\Mu^{-1}(p)$; the point $p$ is called a {\it point
of complete collision}. Points of complete collision lying on edges are
called simply {\it points of collision}. A comotion is said to be {\it
regular} if all the mappings $\alpha_i$ are nondecreasing.

Note that, for a continuous nondecreasing mapping from one circle to
another, the number of simply connected components of the preimage of a
point is always finite, does not depend on the point, and coincides with
the degree of the mapping.

The main property of comotions used in this paper is as follows.

\Lemma 11.
For a regular comotion $\{\alpha_i\}$ on a surface $S$,
$$
(the\ number\ of\ points\ of\ complete\ collision)+\sum_i (1-\deg\alpha_i)
\ge e(S).
$$

In fact, this inequality is a simple corollary of an equality, to
write which we need some additional denotations.

Take a one-to-one continuous orientation preserving mapping
$f\:[0;T)\to\bf T$ from a
half-open interval $[0;T)$ ($T\in\R$) onto a circle.
Let us define the functions $\chi\:\bf T\times\bf T\to\Z$ and
$\psi\:\coprod\limits_{k\in\N} {\bf T}^k\to\Z$ by the formulae
$$
\eqalign{
\chi(x,y)&=\cases{0&if $f^{-1}(x)\le f^{-1}(y)$\cr
                  1&if $f^{-1}(x)>f^{-1}(y)$\cr};\cr
\psi(t_1,\dots,t_k)&=\chi(t_1,t_2)+\chi(t_2,t_3)+\dots+\chi(t_k,t_1);\ \
\psi(t_1)=0.
}
$$

\Lemma 12. The function $\psi$ has the following
properties:
\item{\rm1)}
  it does not depend on the choice of the function $f$ and can be defined
  invariantly as follows: consider an oriented circle $X$, points
  $x_1,\dots,x_k\in X$ lying on this circle in the specified order, and
  a continuous nondecreasing mapping $F\:X\to\bf T$, which maps each point
  $x_i$ to a point $t_i$ and the arc $[x_i,x_{i+1}]$ onto the arc
  $[t_i,t_{i+1}]$ \(here the subscripts are modulo $k$ and an arc with
 coinciding endpoints is assumed to be the singleton\); then
  $\psi(t_1,\dots,t_k)=\deg F$;
\item{\rm2)}
  it takes integer nonnegative values;
\item{\rm3)}
  $\psi(t_1,\dots,t_k)=0$ if and only if all $t_i$ are equal.

\Proof
Properties 2 and 3 are obvious. To prove property 1, note that, by
definition, $\psi(t_1,\dots,t_k)$ equals to the number of the half-open
arcs $(t_i,t_{i+1}]$ containing $f(0)$ (here the subscripts are modulo
$k$ and a half-open arc with coinciding endpoints is assumed to be empty).
This number, in its turn, is the number of $F$-preimages of the point
$f(T-\epsilon)$, where $\epsilon$ is a sufficiently small positive real
number. The lemma is proven.

\Lemma 13.
Let $\Alpha\:\coprod\d D_i\to\bf T$ be a regular comotion on a surface
$S$. Let us define the weight $\nu$ for a face $D_i$, an edge $e$, and
a vertex $v$ by the following formulae:
$$
\eqalign{
\nu(D_i)=&1-\deg\alpha_i,\cr \nu(e)=&-1+\hbox{\rm(the number of connected
components of the set of points of $e$ not being points of collision)},\cr
\nu(v)=&1-\psi(\Alpha(c_0),\dots,\Alpha(c_k)), }
$$
where
$c_0,\dots,c_k$ are all corners at the vertex $v$
enumerated anticlockwise. Then the sum of the weights of all faces, edges,
and vertices equals the Euler characteristic of the surface $S$.
\rm
(Recall that, according to the definition, edge's endpoints do not belong
to the edge.)

\Proof
Note that the total weight of the faces, edges, and vertices is
invariant with respect to subdivision of edges, i.e.,
it remains the same when a new vertex $v$ dividing an edge $e$ into two
edges $e_1$ and $e_2$ is added. Indeed, if the new vertex $v$ is a point
of the edge $e$ at which no collision occurs, then the total number of
connected components of the set of points of the edges $e_1$ and $e_2$ not
being points of collision is large by one than the similar value for the
edge $e$, and the weight of the vertex $v$ is 0. Thus, the total weight
does not change:  $\nu(e_1)+\nu(e_2)+\nu(v)=\nu(e)$. If the new vertex $v$
is a collision point of the edge $e$, then the total number of connected
components of the set of points of the edges $e_1$ and $e_2$ not being
points of collision equals the similar value for the edge $e$, and
the weight of the vertex $v$ is 1. Thus, the total weight does not
change either.

Subdividing the edges, we can achieve the situation when each edge $e$ has
the following properties:
\item{1)}
 at some moment of time, there are no cocars on the closure of $e$;
\item{2)}
 either there occur no collisions on $e$,
 or collisions occur at each point of $e$ (in the latter case, the
 function $\Alpha$ is constant on $\Mu^{-1}(e)$).

Note that the weight of an edge $e$ having these properties equals
either 0 or $-1$ and can be written in the form
$$
\nu(e)=\psi(\Alpha(c_1),\Alpha(c_2),\Alpha(c_3),\Alpha(c_4))-1,
$$
where $c_1,c_2,c_3,c_4$ are the corners adjacent to $e$ enumerated
clockwise (as in Figure 7). The validity of this formula can be easily
verified by using the properties of function $\psi$ (Lemma 12) and the
fact that the properties 1 and 2 of the edge $e$ imply that the intervals
$(\Alpha(c_1),\Alpha(c_2))$ and $(\Alpha(c_3),\Alpha(c_4))$ are disjoint.

According to Lemma 12, the weight of a face can be written in the form
$$
\nu(D_i)=1-\psi(\alpha_i(c_0),\dots,\alpha_i(c_k)),
$$
where $c_0,\dots,c_k$ are all corners of the face $D_i$ enumerated
anticlockwise. To complete the proof, we use the following obvious fact.

\Lemma 14.
Suppose that we have a map on a surface $S$ and two functions $g$
and $h$ assigning numbers to pairs of corners. Let us define the weight
$\nu$ for a face $D$, an edge $e$, and a vertex $v$ by the following
formulae:
$$
\eqalign{ \nu(D)=&1-\sum_{i=0}^kg(c_i,c_{i+1}),\cr
\nu(e)=&-1+g(c_1,c_2)+h(c_2,c_3)+g(c_3,c_4)+h(c_4,c_1),\cr
\nu(v)=&1-\sum_{i=0}^kh(c_i,c_{i+1});\cr
}
$$
in the first formula, $c_0,\dots,c_k$ are all corners of the face $D$
enumerated anticlockwise;
in the second formula, $c_1,\dots,c_4$ are all corners adjacent to the
edge $e$ enumerated clockwise (Fig. 7); and, in the last formula,
$c_0,\dots,c_k$ are all corners at the vertex $v$ enumerated
anticlockwise. Then the total weight of all faces, edges, and vertices
equals the Euler characteristic of the surface $S$.

To prove this fact, it is sufficient to note that, when the weights are
summed, all the values of the functions $g$ and $h$ are cancelled: each
term $g(x,y)$ which occurs in the weight of a face with a minus sign
occurs also in the weight of an edge with a plus sign; similarly, each
term $h(x,y)$ occuring in the weight of a vertex with a minus sign occurs
also in the weight of an edge with a plus sign.

To complete the proof of Lemma 13 we should put
$g(c_1,c_2)=h(c_1,c_2)=\chi(\Alpha(c_1),\Alpha(c_2))$ in Lemma 14.

Lemma 11 immediately follows from Lemma 13; we should only note that
the weight of an edge (from Lemma 13) is not higher than the number of
collision points of this edge, and the weight of a vertex is 1 if
a complete collision occurs at this vertex and is nonpositive otherwise.

\Def
A {\it multiple motion of period $T\in\R$} on a surface $S$
is a map $\{\mu_i\:D_i\to S\ |\ i\in I\}$ on this surface
and a set of mappings
$\{\alpha_{i,j}\:\R\to\d D_i;\ i\in I,j=1,\dots,d_i\}$ (called
cars) satisfying the following periodicity conditions:
\item{1)}
  $\alpha_{i,j}(t+T)=\alpha_{i,j+1}(t)$ for any $t\in \R$ and
  $j\in\{1,\dots,d_i\}$ (subscripts are modulo $d_i$);
\item{2)}
  there exists such a partitioning of
  each circle $\d D_i$ into $d_i$ arcs that during the time interval
  $[0,T]$ each car $\alpha_{i,j}$ moves along the $j$th arc.

\noindent
The positive integers $d_i$ are called the {\it multiplicities} of the
multiple motion. A multiple motion is called {\it regular} if all the
function $\alpha_{i,j}$ are orientation-preserving coverings.

\Example 2.
Let us specify a multiple motion on the map shown in Figure 1 as follows:
the cars $\alpha$, $\beta$, $\gamma$, $\delta$, and $\epsilon$ going
around the faces $A$, $B$, $C$, $D$, and $E$, respectively,
move as in Example 1, i.e., with unit speed (one edge per unit time)
in the positive direction being at time zero at the corners $a_0$,
$b_0$, $c_0$, $d_0$, and $e_0$, respectively; in addition, there is also
a car $\beta'$ going around the face $B$ in the positive direction with unit
speed, being at time zero at the corner $b_3$. In Figure 1, the
positions of the cars at time $t=4/3$ are shown.  This
multiple motion is regular, it has period 3 and set of multiplicities
$\{1,2,1,1,1\}$.  The complete collisions occur at the same three points
as in Example 1 (they are marked by the exclamation signs in Figure 1).
However, in this example, no variation of the motion schedule can decrease
the number of points of complete collision.

\Lemma 15.
The number of points of complete collision of a regular multiple motion
with multiplicities $\{d_i;\ i\in I\}$ on a surface $S$ cannot be smaller
than
$$
e(S)+\sum_{i\in I}(d_i-1).
$$

\Proof
A regular multiple motion $\{\alpha_{i,j}\}$ induces
a comotion $\{\beta_i\:\d D_i\to \R/T\Z\}$. It is
sufficient to put $\beta_i(x)$ to be such $t$ that
$\alpha_{i,j}(t)=x$; this function does not depend on $j$, is
well-defined modulo $T$, and is a covering of degree $d_i$. To complete
the proof, it remains to note that the points of complete collision of
this comotion coincide with the points of complete collision of the
initial multiple motion and apply Lemma 11.

A multiple motion $\{\alpha_{i,j}\}$ is called a {\it multiple motion with
separated stops} if all the cars $\alpha_{i,j}$ are nondecreasing functions
and the stops are separated in the sense of the definition from Section 5.

\Lemma 16.
The number of points of complete collision of a multiple motion with
separated stops on a surface $S$ cannot be smaller than
$$
e(S)+\sum_{i\in I}(d_i-1),
$$
where $d_i$ are the multiplicities of the multiple motion.

This lemma is derived from Lemma 15 in exactly the same way as Lemma 4 is
derived from Lemma 3.

We say that a map on a surface is {\it $2$-graded} if it has the following
additional structure: some vertices and some faces are distinguished and
called {\it exterior}, the remaining vertices and faces are called {\it
interior}; in addition, some faces are called {\it large}, the remaining
faces are called {\it small}.

Suppose that we have a (multiple) motion on a 2-graded map $\Mu$ on a
surface $S$.  We say that two different large faces $A$ and $B$ {\it badly
contact each other} if they have such corners $a_1,a_2\in\d A$ and
$b_1,b_2\in\d B$ that
\item{1)}
 $a_i$ and $b_i$ are nonadjacent corners at some interior vertex
 of complete collision $v_i=\Mu(a_i)=\Mu(b_i)$ (for $i=1,2$);
\item{2)}
 $v_1\ne v_2$;
\item{3)}
 the closed path on the surface $S$ formed by the fragment
 $[a_1,a_2]$ of the boundary of the face $A$ and the fragment $[b_2,b_1]$
 of the boundary of the face $B$ does not meet exterior vertices and
 bounds a disk submap all of whose vertices are interior and all
 cells are small and interior.

\noindent
We say that a large face $A$ {\it badly contacts itself} if either
it has such four different corners $a_1,a_2,b_2,b_1\in\d A$ (disposed in
this order) that the above conditions 1), 2), and 3) hold, or it has
corners $a,b\in\d A$ being different nonadjacent corners at some interior
vertex of complete collision $v=\Mu(a)=\Mu(b)$ and the fragment $[a,b]$ of
the boundary of $A$ does not meet exterior vertices on the surface $S$
and bounds a disk submap all of whose vertices are interior and all
cells are small and interior.

Figure 9 illustrates all cases of bad contact (contiguity); the asterisks
mark submaps containing no exterior vertices, no exterior
faces, and no large faces; at the vertices marked by the exclamation
signs, complete collisions occur.

\input 9.PIC

\Lemma 17. No spherical 2-graded map with a single
exterior vertex and without exterior faces admits a multiple
motion with separated stops satisfying the following conditions:
\item{\rm1)}
  the multiplicity of the motion is at
  least 4 on each large face and at least 1 on each small face;
\item{\rm2)}
  at each interior point of complete collision, there are at least 2
 nonadjacent corners of large cells;
\item{\rm3)}
 there are no badly contacting (themselves or others) large cells.

\Proof
To obtain a contradiction with Lemma 16, it is sufficient to show that
conditions 2) and 3) imply that
$$
(the\ number\ of\ interior\ points\ of\ complete\ collision)<
3\cdot(the\ number\ of\ large\ faces)+1.
$$
Let us take a point $v_D$ inside each large face $D$. At each interior
point of complete collision $p$, let us choose 2 nonadjacent corners $a$
and $b$ lying on large faces, which we call $A$ and $B$, respectively.
Let us draw an arc on the sphere from the point $v_A$ through the interior
of the face $A$ via the corner $a$ to the point $p$ and, further,
via the corner $b$ through the interior of the face $B$ to the point
$v_B$. If we draw these arcs in such a way that they do not intersect each
other (this is possible, of course), then we obtain a graph $\Gamma$ on
the sphere. The number of vertices of this graph equals the number of
large cells of the map, and the number of edges of the graph equals the
number of interior points of complete collision. To complete the proof, it
remains to note that condition 3) means that the graph $\Gamma$ satisfies
the conditions of the following lemma.

\Lemma 18.
Suppose that we have a finite \(not necessarily connected\) graph
$\Gamma$ on the sphere $S^2$ such that the perimeter of each simply
connected component of $S^2\setminus\Gamma$, except maybe one,
is at least 3. Then the number of edges of this graph is not larger than
the tripled number of its vertices.\hfil\break
\rm Here, the perimeter of a region is the number of edges on the boundary
of this region, where an edge is counted twice if the region lies on
both sides from this edge.

\Proof
If the graph has no edges, then we have nothing to prove.
Assuming that edges exist, let us use induction on the number of connected
components of $\Gamma$. If the graph is connected, then the Euler formula
gives $V-E+F=2$, where $V$, $E$, and $F$ are the numbers of vertices, edges,
and faces of the corresponding map on the sphere. Adding one
vertex inside the exceptional component and joining it by a new edge with
a vertex on the boundary of this component, we obtain a map each of whose
face is at least a triangle. Therefore, $F\le {2\over3}(E+1)$.  Hence,
$V-E+{2\over3}E=V-{1\over3}E\ge2-{2\over3}={4\over3}$. Thus, for
connected graphs, the lemma is true.

If edges exist, but the graph is not connected, then we join a connected
nonsingleton component with another component by a new edge.
We thereby decrease the number of connected components without changing the
number of vertices, and increase the number of edges; the new simply
connected component of the complement to the graph has perimeter at least
three. The application of the induction hypothesis completes the proof of
Lemma 18 and also of Lemma 17.

We say that a map with oriented edges is of type $B_m$, $m\ge0$,
if the sequence of pre-edge orientations of each face has one of the
following three forms:
\item{a)} $+(+-)^{m+1}$ (Figure 6b);
\item{b)} $-(-+)^{m+1}$ (Figure 6c);
\item{c)} $((+)^{k+1}(-)^{l+1})^s$, where $k,l,s\ge1$ (Figure 10).

\input 10.PIC

Let us define the {\it standard multiple motion} on a map of type
$B_m$ as follows. On the faces of types $+(+-)^{m+1}$ and $-(-+)^{m+1}$,
the standard motion (of multiplicity 1) coincides with the standard
motion on such faces from the definition of the standard motions on
maps of type $A_m$ (Fig. 6b and 6c). On a face of type
$((+)^{k+1}(-)^{l+1})^s$, we define an $s$-multiple motion as the lift of
the standard motion on a face of type $(+)^{k+1}(-)^{l+1}$ (from the
definition of the standard motion on a map of type $A_m$). To be more
precise, there is a natural pre-edge orientation preserving
$s$-fold covering $\pi\:\d B\to\d A$, where $B$ is a face of type
$((+)^{k+1}(-)^{l+1})^s$ and $A$ is a face of type
$(+)^{k+1}(-)^{l+1}$; the standard car $\alpha\:\R\to\d A$ has precisely
$s$ liftings, i.e., such functions $\alpha_1,\dots,\alpha_s\:\R\to\d B$
that $\alpha=\pi\alpha_i$; these functions are called the standard cars
moving around $B$.

Figure 10 presents the schedule of the standard multiple motion on a face
of type $((+)^{k+1}(-)^{l+1})^s$ for $k=l=2$ and $s=4$: the
number near a vertex is the time ($\mod 4m+2$) when one of the four cars
visits this vertex. In this schedule, we assume, of course, that $m>0$;
for $m=0$, all the 4 cars move uniformly with speed $k+1=l+1=3$ and visit
the vertices marked by the number 0 at time zero.

\Lemma 19. The standard motion on a map of type $B_m$ is
a multiple motion with separated stops. Complete collisions can occur only
at vertices being sources or sinks and only at integer moments of time.

The proof of this lemma is a word-by-word repetition of the proof of Lemma 9.


\s 10. Proof of the main theorem. The difficult case

As above, we assume that the group $G$ is torsion-free and conditions (1)
hold. It is sufficient to prove that $\~G$ is nonsimple.
The group $\~G$ has presentation (2), where we assume $P=\1$,
i.e., presentation (2) has the form
$$
\~G\iso\gp{G, t\ \Biggm|\
ct\prod_{i=0}^m(b_ia_i^t)=1},\quad m\ge0.
\eqno{(4)}
$$

\Lemma 20. $G^t\cap G=\1$ in $\~G$.

\Proof
Assuming the contrary, by Lemma 5 (and the remark after Lemma 10) we
conclude that there exists a reduced spherical diagram over presentation
(4) with a single exterior face, and the label of this face is $g^th$,
where $g,h\in G\setminus\1$.  This diagram is a map of type $A_m$. It has
one (exterior) cell of type $+-$ (Fig.6a), the remaining cells are of
types $+(+-)^{m+1}$ (Fig. 6b) and $-(-+)^{m+1}$ (Fig. 6c). Consider
the standard motion of period $4m+2$ on this map. According to the remark
after Lemma 10, complete collisions can occur only at vertices $A$ and
$B$ lying on the boundary of the exterior face. By Lemma 4, complete
collisions must happen at both these vertices. By Lemma 9, collisions can
occur only in sources
and sinks. Therefore, the vertex $A$, at which the $(-+)$-corner of the
exterior face lies, is a source, and the vertex $B$, at which the
$(+-)$-corner of the exterior face lies, is a sink. For $m=0$, this
implies immediately that the labels of all interior corners at the vertex
$A$ are $a_0^{\pm1}$ (because $A$ is a source), but not all of these
labels are equal (because $B$ is a sink) (Fig. 11); hence, the diagram
contains a reducible pair of cells and is not reduced.

\input 11.PIC

If $m>0$, then complete collisions cannot happen at both vertices $A$ and
$B$. Indeed, suppose that at the vertex $A$ a complete collision occurs at
time $t$. Then one car, moving around an interior face visits the
vertex $B$, at the moment $t_1=t+1$, and another car visits $B$ at the
moment $t_2=t-1$. But $t_1$ never equals $t_2$ modulo the motion period,
because, for $m>0$, the period is larger than 2. Thus, no complete
collision occurs at $B$. This contradiction completes the proof.

\Lemma 21. For any $g,h\in G$, we have
{either} $\gp{g}^{t^2}\cap\gp{h}=\1$ in $\~G$
{or} $\gp{g}^{t^3}\cap\gp{h}=\1$ in $\~G$.

\Proof
Assume the contrary:
$$
g^{kt^2}=h^l,\quad g^{k't^3}=h^{l'}.
$$
Let us raise these equalities to the powers $k'$ and $k$, respectively:
$$
g^{kk't^2}=h^{lk'},\quad g^{kk't^3}=h^{kl'}.
$$
Conjugating the first equality by $t$, we obtain
$$
h^{lk't}=h^{kl'},
$$
whence $h=1$ by virtue of Lemma 20 and the absence of torsion in $G$.

\Lemma 22.
There exists such $d\in\{2,3\}$ that
$$
u\equiv\prod_{i=1}^s y_ix_i^{t^d}\ne1\ {\rm in}\ \~G
$$
for any positive integer $s$ and any $x_i, y_i\in G$
such that
$\Bigl|\{i\ |\ x_i\in\gp{a_m}\}\Bigr|+
\Bigl|\{i\ |\ y_i\in\gp{b_0}\}\Bigr|\le2$ and $u\ne1$ in $G*\gp{t}_\infty$.

\Proof
Let $d$ be such number that $\gp{a_m}^{t^d}\cap\gp{b_0}=\1$; such a
$d\in\{2,3\}$ exists by Lemma 21.

Proving by contradiction, consider a counterexample with the minimal
possible $s$.
By Lemma 5 (and the remark after Lemma 10), there exists a diagram over
presentation (4) on the sphere with a single exterior face whose label
is~$u$.

First, let us show that
{at each vertex there is at most one corner of type $(+-)$ of the
exterior face.} Indeed, supposing that exterior face's corners labelled
by $y_1$ and $y_r$ are corners at the same vertex and considering the
corresponding subdiagrams, we obtain that the equality $u=1$ decomposes
into a product of two equalities
$$
\left(\prod_{i=1}^{r-1}
y_ix_i^{t^d}\right)g=1\quad\hbox{and}\quad g^{-1}\prod_{i=r}^{s}
y_ix_i^{t^d}=1,\quad\hbox{where $g\in G$},
$$
(see Fig. 12, where $s=5$, $d=2$, and $r=3$), at least one of which
contradicts the minimality of the counterexample. Similarly, we can show
that
{at each vertex, there is at most one corner of type $(-+)$ of the
exterior face}.

\input 12.PIC

Note that the diagram is a map of type $B_m$. Consider the standard
multiple motion on this map. According to the remark after Lemma 10,
complete collisions can occur only at vertices lying on the
boundary of the exterior face. Moreover, by Lemma 19, a vertex of complete
collision must be a source or a sink. Therefore, at each vertex of
complete collision, there is exactly one corner of the exterior face, and
a complete collision can occur only at the moments 0 and $2m+1$ ($\mod
4m+2$), because only at these moments the cars moving around the exterior
face visit corners of types $(+-)$ and $(-+)$. Since in these moments all
the remaining cars are at corners labelled by $b_0^{\pm1}$ (at the moment
0) and $a_m^{\pm1}$ (at the moment $2m+1$), the label of exterior face's
corner at a vertex of complete collision must belong to $\gp{b_0}$ (for
a corner of type $(+-)$) or to $\gp{a_m}$ (for a corner of type
$(-+)$). But the exterior face has at most two such corners. Hence, there is
at most two points of complete collision. But according to Lemma 16, there
must be at least $s+1$ such points. Therefore,
$s=1$, and $u$ has the form $b_0^ka_m^{lt^d}$; but such a word does not
equal 1 in $\~G$ by virtue of the choice of $d$.  The lemma is proven.

\Lemma 23. There exists such $d\in\{2,3\}$ that
the presentation
$$ \~G\iso\gp{G, t\ \Biggm|\ ct\prod_{i=0}^m(b_ia_i^t)=1,\
(a^{t^d}b)^4=1},\quad\hbox{where $m\ge0$},
$$
is aspherical for any elements $a,b\in G$ such that $a^2\notin\gp{a_m}$
and $b^2\notin\gp{b_0}$.

\Proof
Let $d$ be the number, whose existence is asserted by Lemma 22.

Proving by contradiction, consider a spherical reduced diagram over this
presentation with a single exterior vertex and without exterior faces.
This diagram is a map of type $B_m$. Consider the standard multiple
motion on this map.  Let us call the cells with boundary label
$(a^{t^d}b)^{\pm4}$ large, the other cells will be called small.

Let us show that conditions 1), 2), and 3) from Lemma 17 hold.

Condition 1) holds by the definition of the standard multiple motion.

Let us show that condition 2) holds. At each vertex of complete collision,
there is at least one corner of a large cell (see the remark after Lemma 10).
By Lemma 19, a vertex of complete collision must be a source or a sink;
such vertices are visited by cars moving around large cells only at
the moments 0 and $2m+1$ ($\mod 4m+2$), because only at these moments
a car moving around a large cell visits corners of types $(+-)$ and
$(-+)$. At these moments, the cars moving around small cells are at
corners labelled $b_0^{\pm1}$ (at the moment 0) and $a_m^{\pm1}$ (in
the moment $2m+1$). Therefore, the label of each corner at a vertex of
complete collision is $a_m^{\pm1}$ or $a^{\pm1}$ if this
vertex is a source and $b_0^{\pm1}$ or $b^{\pm1}$ if this
vertex is a sink. The labels of adjacent corners are not mutually
inverse, because the diagram is reduced. Suppose that condition 2) does
not hold, i.e., at some interior vertex of complete collision, there is no
pair of different nonadjacent corners of large cells, that is, corners with
labels $a^{\pm1}$ or $b^{\pm1}$. Suppose that this vertex is a source
(the cast of a sink is quite similar).
Then its label has the form
$$
\hbox{either } a^{\pm1}a_m^k
\hbox{\quad or } a^{\pm2}a_m^k,
\hbox{\quad or } a^{\pm3}.
$$
On the other hand, this label must be 1 in $G$, because it is the label
of an interior vertex. Thus, we have three equalities, one of which must
hold in the group $G$. The first two of these equalities contradict the
condition $a^2\notin\gp{a_m}$, and the third contradicts the absence of
torsion in $G$.
This shows that condition 2) of Lemma 17 holds.

Suppose that condition 3) does not hold. We have a disk subdiagram
(not containing the exterior vertex)
whose interior cells are small and the label of the contour has the form
$$
u=\prod_{i=1}^s y_ix_i^{t^d},
$$
where $s\in\N$, $x_i, y_i\in G\setminus\1$,
$\bigl|\{i\ |\ x_i\ne a^{\pm1}\}\bigr|+
\bigl|\{i\ |\ y_i\ne b^{\pm1}\}\bigr|\le2$, and one or two exceptional
coefficients are nonzero powers of $b_0$ or $a_m$ (Fig. 13). This
means (by Lemma 5) that $u=1$ in $\~G$, which is impossible by Lemma 22.

\input 13.PIC

Thus, the diagram under consideration has properties 1), 2), and 3) from
Lemma 17, which asserts that such diagrams do not exist. This
contradiction completes the proof of Lemma 23.

The assertion of the main theorem (in the difficult case) follows from
Lemmata 23 and 7 and the obvious fact that if the squares of all elements
of the nontrivial group $G$ lie in the same cyclic subgroup, then $G$ does
not coincide with its commutator subgroup (because it is metabelian) and,
therefore, $\~G$ does not coincide with its commutator subgroup either.

\s{\rm REFERENCES}

\item{[CG95]}
Clifford A., Goldstein R. Z.
Tesselations of $S^2$ and equations over torsion-free groups
{// Proc. Edinburgh Math. Soc.} 1995. {V.38}. P.485--493.

\item{[CG00]}
Clifford A., Goldstein R. Z.
Equations with torsion-free coefficients
{// Proc. Edinburgh Math. Soc.} 2000. {V.43}. P.295--307.

\item{[CR01]}
Cohen M. M., Rourke C.
The surjectivity problem for one-generator, one-relator extensions of
torsion-free groups
{// Geometry \& Topology}. 2001. {V.5}. P.127--142.

\item{[FeR96]}
Fenn R., Rourke C.
Klyachko's methods and the solution of equations over torsion-free groups
{// L'Enseignment Math\'ematique.} 1996. {T.42}. P.49--74.

\item{[FeR98]}
Fenn R., Rourke C.
Characterisation of a class of equations with solution over torsion-free
groups,
from {``The Epstein Birthday Schrift"},
{(I. Rivin, C. Rourke and C. Series, editors)},
{Geometry and Topology Monographs.} 1998. {V.1}. P.159-166.

\item{[FoR03]}
Forester M., Rourke C.
Diagrams and the second homotopy group.
{// arXiv:math.AT/0306088}. v.1. 5 Jun 2003

\item{[GR62]}
Gerstenhaber M., Rothaus O. S.
The solution of sets of equations in groups
{//  Proc. Nat. Acad. Sci. USA}. 1962. {V.48} P.1531--1533.

\item{[H83]}
Howie J.
The solution of length three equations over groups
{// Proc. Edinburgh Math. Soc.} 1983. {V.26}. P.89--96.

\item{[K93]}
Klyachko Ant. A.
A funny property of a sphere and equations over groups
{// Comm. Algebra}. 1993. {V.21}. P.2555--2575.

\item{[K94]}
Klyachko Ant. A.
Kervaire--Laudenbach conjecture and equation over groups,
{Cand. Sci. Dissertation},
{Moscow: MSU,} 1994.

\item{[K97]}
Klyachko Ant. A.
Asphericity tests
{// J. Algebra}. 1997. {V.7}. P.415--431.

\item{[KP95]}
Klyachko Ant. A., Prishchepov M. I.
The descent method for equations over  groups
{// Moscow Univ. Math. Bull.} 1995, {V.50}  P. 56--58.

\item{[LS77]}
Lyndon R. C., Schupp P. E.
{Combinatorial group theory},
Springer-Verlag, Berlin/Heidelberg/New~York, 1977.

\item{[MKS66]}
Magnus W., Karrass A., Solitar D.
{Combinatorial group theory},
Interscience Publishers, John Wiley \& Sons, New~York/London/Sidney, 1966.

\end